# Stability of Triangular Equilibrium Points in Robe's Generalised Restricted Three Body Problem

**K.T.Singh(i), B.S.Kushvah(ii) and B.Ishwar(iii)**

(i) Lecturer in Mathematics M.B.College ,Imphal(Manipur),
(ii) J.R.F D.S.T. Project*,
(iii)Principal Investigator D.S.T.Project*,

*University Department of Mathematics ,
B.R.A. Bihar University Muza_arpur-842001,
Email:ishwar bhola@hotmail.com

**Abstract :** We have examined the stability of triangular equilibrium points in Robe's generalized restricted three body problem .The problem is generalised in the sense that more massive primary has been taken as an oblate spheroid.We have found the position of triangular equilibrium points .We have obtained variational equations of the problem.With the help of characteristic roots, we conclude that triangular points are unstable .Robes's result may be verified from the generalised result.

**Introduction :**

This paper is devoted to the study of location of triangular equilibrium points and its linear stability in the Robe's restricted three body problem with perturbation coriolis and centrifugal forces and the mass $m_1$ is taken as an oblate spheroid. Robe's Restricted Three Body Problem, is new kind of restricted three body problem in which one of the primaries is a rigid spherical shell $m_1$ filled with a homogeneous incompressible fluid of density $r_1$. The second primary is a mass point $m_2$ outside the shell and the third body $m_3$ is a small solid sphere of density $r_3$ inside the shell with the assumption that the mass and radius of $m_3$ describes a Kaplerian orbit around it. Further, he has discussed the linear stability of the equilibrium point for the whole range of parameters occurring in the problem. Later on , Shrivastava and Garain (1991) have studied the effect of small perturbations in the coriolis and centrifugal forces on the location of equilibrium points in the Robe's problem. Plastino and Plastino (1995) have considered the Robes problem by taking the shape of the fluid body as Roche's ellipsoid. They have studies the linear stability of the equilibrium solution too. Giordano, Plastino and Plastino (1997) have discussed the effect of drag force on the existence and stability of the equilibrium points in Robe's problem. Hallan and Rana established that in the Robe's elliptic restricted three body problem there is only one equilibrium point while in circular problem there are two, three and infinite number of equilibrium points.

In this paper we have found the position of triangular equilibrium points. They lie in xz – plane. Also we have examined the stability of triangular equilibrium points. We have found the variational equations. The values of second order partial derivatives were found at equilibrium points. Then the characteristic equation was obtained. Finally, we conclude that triangular equilibrium points are unstable.

**LOCATION OF TRIANGULAR EQUILIBRIUM POINT :**

The equation of motion are

$$\left.\begin{array}{l}\Omega_x = \ddot{x} - 2n\dot{y} \\ \Omega_y = \ddot{y} + 2n\dot{x} \\ \Omega_z = \ddot{z}\end{array}\right\} \quad \text{---- (1)}$$

where $\Omega = \frac{n^2}{2}\left(x^2+y^2\right) - kr_1^2 + \frac{m}{r_2}$ , $r_1^2 = (x+m)^2 + y^2 + z^2$ ,

$r_2^2 = (x+m-1)^2 + y^2 + z^2$ , $k = \frac{4}{3}p\,r_1\left(1 - \frac{r_1}{r_3}\right)$, $m = \frac{m_2}{m_1+m_2}$ , $n^2 = 1 + \frac{3}{2}A_1$

Co – ordinates of $m_1$ and $m_2$ are $(-m,0,0)$ and $(1-m,0,0)$. The libration points exist when

$\Omega_x = \Omega_y = \Omega_z = 0$ ie.

$$\left.\begin{array}{l} n^2x - 2k(x+m) - \frac{m}{r_2^3}(x+m-1) = 0 \\ n^2y - 2ky - \frac{m}{r_2^3}y = 0 \\ -2kz - \frac{mz}{r_2^3} = 0 \end{array}\right\} \quad \text{---- (2)}$$

from second equation we get $\left(n^2 - 2k - \frac{m}{r_2^3}\right)y = 0$ ,

$y = 0$ $\quad \because n^2 - 2k - \frac{m}{r_2^3} \neq 0$. From third equation of equation (2) we get,

$$-2k = \frac{m}{r_2^3} \quad [\because z \neq 0] \quad \text{---- (3)}$$

From first equations of equations (2) and equations (3) we get,

or $n^2x - 2k(x+m) + 2k(x+m-1) = 0$

or $x = \frac{2k}{n^2}$ we know $r_2^2 = (x+m-1)^2 + y^2 + z^2 = (x+m-1)^2 + z^2 \quad [\because y = 0]$

Putting the values of x and $r_2 = \left(-\frac{m}{2k}\right)^{1/3}$, we get $\left(-\frac{m}{2k}\right)^{2/3} = \left(\frac{2k}{n^2} + m - 1\right)^2 + z^2$

or $\quad z^2 = -\left(\frac{2k}{n^2} + m - 1\right)^2 + \left(-\frac{m}{2k}\right)^{2/3}$

or $$= \pm\left[\left(-z\frac{m}{2k}\right)^{2/3} - \left(\frac{2k}{n^2} + m - 1\right)^2\right]^{1/2}$$ ---- (4)

or $$z = \pm\left(b_1^2 - a_1^2\right)^{1/2}$$ ---- (5)

where $$b_1 = \left(-\frac{m}{2k}\right)^{1/3}, \qquad a_1 = \frac{2k}{n^2} + m - 1$$

Hence co -ordinates of triangular equilibrium points are $(x,o,z) = \left(\frac{2k}{n^2}, 0, \pm\left(b_1^2 - a_1^2\right)^{1/2}\right)$, provided k < 0, and $\frac{2k}{n^2} + m > 0$ , that is in the $m - \frac{2k}{n^2}$ plane the point $\left(m, \frac{2k}{n^2}\right)$ lies within the triangular region bounded by the lines k = 0, μ = 1, $\frac{2k}{n^2} + m = 0$.

**STABILITY OF THE TRIANGULAR POINTS :**

Equations of motion are

$$\left.\begin{aligned} \ddot{x} - 2n\dot{y} &= \frac{\partial \Omega}{\partial x} \\ \ddot{y} + 2n\dot{x} &= \frac{\partial \Omega}{\partial y} \\ \ddot{z} &= \frac{\partial \Omega}{\partial z} \end{aligned}\right\}$$ --- (6)

where

$$\Omega = \frac{n^2}{z}\left(x^2 + y^2\right) - kr_1^2 + \frac{m}{r_2}$$ ---- (7)

$$\left.\begin{aligned} r_1^2 &= (x+m)^2 + y^2 + z^2 \\ r_2^2 &= (x+m-1)^2 + y^2 + z^2 \\ k &= \frac{4}{3}p r_1\left(1 - \frac{r_1}{r_3}\right) \end{aligned}\right\}$$ ---- (8)

Triangular equilibrium points are given by $\left(\frac{2k}{n^2}, 0, \pm\left(b_1^2 - a_1^2\right)^{1/2}\right)$

where $a_1 = \frac{2k}{n^2} + m - 1,$ $b_1 = \left(-\frac{m}{2k}\right)^{1/3}$

Let $x, h, z$ be the small displacements in the equilibrium point $(x_0, y_0, z_0)$. Then the third body will be displaced to $(x_0 + x, y_0 + h, z_0 + z)$ and the variational equations in linearised form are

$$\left.\begin{aligned} \ddot{x} - 2n\dot{h} &= \Omega^0_{xx} x + \Omega^0_{xy} h + \Omega^0_{xz} z \\ \ddot{h} + 2n\dot{x} &= \Omega^0_{yx} x + \Omega^0_{yy} h + \Omega^0_{yz} z \\ \ddot{z} &= \Omega^0_{zx} x + \Omega^0_{zy} h + \Omega^0_{zz} z \end{aligned}\right\} \quad \text{---- (9)}$$

we have, $\Omega = \frac{n^2}{2}(x^2 + y^2) - kr_1^2 + \frac{m}{r_2}$, At the triangular equilibrium points $(x, o, z) = \left(\frac{2k}{n^2}, o, \pm(b_1^2 - a_1^2)^{1/2}\right)$

We have, $\Omega^0_{xx} = \left[n^2 - 2k - \frac{m}{r_2^3} + \frac{3m}{r_2^5}(x + m - 1)^2\right],$

$$= n^2 - 6k\left(\frac{a_1}{b_1}\right)^2, \ \Omega^0_{xx} = n^2 - 6k\left(\frac{a_1}{b_1}\right)^2$$

$$\Omega^0_{xy} = \frac{3m(x + m - 1)y}{r_2^5} \qquad = 3\left(\frac{m}{r_2^3}\right)\frac{(x + m - 1)y}{r_2^2} = 0$$

similarly, $\Omega^0_{zy} = 3\frac{m}{r_2^5} yz = 0$, $\Omega^0_{yz} = 0$, $\Omega^0_{yx} = 0$,

$$\Omega^0_{xz} = \frac{3m(x + m - 1)z}{r_2^5} \qquad = -6k\frac{q_1(b_1^2 - q_1^2)^{1/2}}{b_1^2} = \Omega^0_{zx}$$

$$\Omega^0_{yy} = \left[n^2 - 2k - \frac{m}{r_2^3} + \frac{3my^2}{r_2^5}\right] = n^2 - 2k - (-2k) + 0$$

$$\Omega^0_{yy} = n^2 \qquad , \Omega^0_{zz} = -2k - \frac{m}{r_2^3} + \frac{3mz^2}{r_2^5} \ , \Omega^0_{zz} = -2k + 2k + 3(-2k)\frac{z^2}{r_2^2}$$

$$= -6k\frac{\left(b_1^2 - a_1^2\right)}{b_1^2},\ \Omega_{zz}^0 = \frac{-6k\left(b_1^2 - a_1^2\right)}{b_1^2}$$

We suppose that $x = Ae^{lt},\ \ h = Be^{lt}\ \ and\ \ z = Ce^{lt}$

Substituting these values in equation (4.3.4) , we get

$$Al^2 - 2nlB = \Omega_{xx}^0 A + \Omega_{xz}^0 C,\quad Bl^2 + 2nlA = \Omega_{yy}^0 B,\quad Cl^2 = \Omega_{zx}^0 A + \Omega_{zz}^0 C$$

or
$$\left(l^2 - \Omega_{xx}^0\right)A - 2nlB - \Omega_{xz}^0 C = 0$$
$$2nlA + \left(l^2 - \Omega_{yy}^0\right)B + 0.C = 0$$
$$-\Omega_{zx}^0 A + 0.B + \left(l^2 - \Omega_{zz}^0\right)C = 0$$

Above system of equation has singular solution if

$$\begin{vmatrix} \left(l^2 - \Omega_{xx}^0\right) & -2nl & -\Omega_{xz}^0 \\ 2nl & \left(l^2 - \Omega_{yy}^0\right) & 0 \\ -\Omega_{zx}^0 & 0 & \left(l^2 - \Omega_{zz}^0\right) \end{vmatrix} = 0$$

ie., $$\left[\left(l^2 - \Omega_{xx}^0\right)\left(l^2 - \Omega_{yy}^0\right) + 4n^2 l^2\right]\left(l^2 - \Omega_{zz}^0\right) - \left(l^2 - \Omega_{yy}^0\right)\left(\Omega_{zx}^0\right)^2 = 0$$

Putting the values of $\Omega_{xx}^0$, $\Omega_{yy}^0$, $\Omega_{zz}^0$ and $\Omega_{zx}^0$ in above expression, we have

$$l^6 + 2\left(n^2 + 3k\right)l^4 + n^2\left[n^2 - 6k\frac{\left(3a_1^2 - b_1^2\right)}{b_1^2}\right]l^2 + 6n^4 k\frac{\left(b_1^2 - a_1^2\right)}{b_1^2} = 0$$

$$^6 + pl^4 + ql^2l + r = 0 \qquad \text{---- (10)}$$

where $p = 2\left(n^2 + 3k\right),\ q = n^2\left[n^2 - 6k\frac{\left(3a_1^2 - b_1^2\right)}{b_1^2}\right],\ r = 6n^4 k\frac{\left(a_1^2 - b_1^2\right)}{b_1^2}$

Since $k < 0\ \ and\ \frac{b_1^2 - a_1^2}{b_1^2} = 1 - \left(\frac{a_1}{b_1}\right)^2 > 0$ ,implies that r is negative and $p > 0$ , $q > 0$ .

Let $f\left(l\right) = l^6 + pl^4 + ql^2 + r$ ---- (11)

Since $f\left(l\right)$ has one change of sign and $f\left(-l\right)$ also has one change of sign so there are two real roots one is positive and other is negative.

$f(l)$ is a polynomial of even degree with real constant term so there are at least two real roots one positive and other negative.Hence triangular points are unstable due to one positive root of equation (11) .Finally we conclude that triangular equilibrium points are unstable.